\documentclass[a4paper,12pt]{amsart}
\usepackage{hyperref,fancyhdr,mathrsfs,amsmath,amscd,amsthm,amsfonts,latexsym,amssymb,stmaryrd}
\usepackage[all]{xy}

\DeclareMathOperator{\grad}{grad}

\DeclareMathOperator{\dif}{d}

\DeclareMathOperator{\Lie}{\mathcal{L}}

\renewcommand{\H}{\mathscr{H}}
\newcommand{\V}{\mathscr{V}}
\newcommand{\E}{\mathscr{E}}
\newcommand{\F}{\mathscr{F}}

\newcommand{\J}{\mathcal{J}}

\def \a{\alpha}

\def \b{\beta}

\def \ep{\varepsilon}

\def \G{\Gamma}

\def \o{\omega}

\def \phi{\varphi}
\def \Phi{\varPhi}
\def \p{\pi}

\def \s{\sigma}

\def \Hq{\mathbb{H}\,}
\def \C{\mathbb{C}\,}

\def\widecheckg{g^{\hspace*{-2.5pt}\vbox to 5pt{\hbox to
0pt{\LARGE$\check{}$}}}\hspace*{2pt}}

\def\widecheckl{\lambda^{\hspace*{-3.5pt}\vbox to 8pt{\hbox to
0pt{\LARGE$\check{}$}}}\hspace*{2pt}}

\begin{document}

\title{On the local structure of generalized K\"ahler manifolds}
\author{Liviu~Ornea and Radu~Pantilie}
\thanks{The authors are  partially supported by a PN II IDEI Grant, code 1193.}
\email{\href{mailto:lornea@gta.math.unibuc.ro}{lornea@gta.math.unibuc.ro},
       \href{mailto:Radu.Pantilie@imar.ro}{Radu.Pantilie@imar.ro}}
\address{L.~Ornea, Universitatea din Bucure\c sti, Facultatea de Matematic\u a, Str.\ Academiei nr.\ 14,
70109, Bucure\c sti, Rom\^ania, \emph{and}
Institutul de Matematic\u a ``Simion~Stoilow'' al Academiei Rom\^ane,
C.P. 1-764, 014700, Bucure\c sti, Rom\^ania}
\address{R.~Pantilie, Institutul de Matematic\u a ``Simion~Stoilow'' al Academiei Rom\^ane,
C.P. 1-764, 014700, Bucure\c sti, Rom\^ania}

\newtheorem{thm}{Theorem}[section]
\newtheorem{lem}[thm]{Lemma}
\newtheorem{cor}[thm]{Corollary}
\newtheorem{prop}[thm]{Proposition}

\theoremstyle{definition}

\newtheorem{defn}[thm]{Definition}
\newtheorem{rem}[thm]{Remark}
\newtheorem{exm}[thm]{Example}

\numberwithin{equation}{section}

\maketitle
\thispagestyle{empty}

\vspace{-7pt}

\section*{Abstract}
\begin{quote}
{\footnotesize  Let $(g,b,J_+,J_-)$ be the bihermitian structure corresponding to a generalized K\"ahler structure.
We find natural integrability conditions, in terms of the eigendistributions of $J_+J_-+J_-J_+$\,, 
under which $\dif\!b=0$\,.}
\end{quote}

\section*{Introduction}

\indent
A \emph{generalized almost complex structure} on a smooth (connected) manifold is given by a vector subbundle
$L\subset\bigl(TM\oplus T^*M\bigr)^{\C}$ such that $L\cap\overline{L}=\{0\}$ and which is maximally isotropic with respect
to the canonical inner product
$$<\!X+\a, Y+\b\!> = \tfrac 12\bigl(\a(Y)+\b(X)\bigr)\;.$$
\indent
If $E=\pi_{TM}(L)$ is a bundle, where $\p_{TM}:TM\oplus T^*M\to TM$ is the projection,
then there exists a unique complex two-form $\ep\in\G(\Lambda^2E^*)$ such that
$L=L(E,\ep)$\,, where
$$L(E,\ep)=\{X+\a\, |\, X\in E, \a|_E=\ep(X)\}\;.$$
Furthermore, by \cite{Gua-thesis}\,, to which we refer for all of the facts on generalized complex structures recalled here,
the condition $L\cap\overline{L}=\{0\}$ is equivalent to $E+\overline{E}=T^{\C}\!M$ and
$\mathrm{Im}\bigl(\ep|_{E\cap \overline{E}}\bigr)$ is non-degenerate.\\
\indent
A generalized almost complex structure $L$ is \emph{integrable} if its space of sections is closed under the \emph{Courant bracket},
defined by
$$[X+\a, Y+\b]=[X, Y] + \Lie_X\!\b - \Lie_Y\!\a - \tfrac 12 \dif(\iota_X \b-\iota_Y\a)\;,$$
for any $X+\a,Y+\b\in\G(L)$\,.\\
\indent
A \emph{generalized complex structure} is an integrable generalized almost complex structure.
Obviously, any generalized complex structure corresponds to a linear complex structure on $TM\oplus T^*M$ whose eigenbundle,
corresponding to ${\rm i}$\,, is isotropic, with respect to the canonical inner product, and its space of sections is closed
under the Courant bracket.\\
\indent
A generalized almost complex structure of the form $L=L(E,\ep)$ is integrable if and only if the space of sections
of $E$ is closed under the (Lie) bracket and $\dif\!\ep(X,Y,Z)=0$\,, for any $X,Y,Z\in E$.\\
\indent
A particular feature of Generalized Complex Geometry is that imposing Hermitian compatibility to a generalized almost complex structure
and a Riemannian metric on $TM\oplus T^*M$, compatible with the canonical inner product, forces the manifold to admit a second generalized
almost complex structure, commuting with the first one. One arrives to the notion of \emph{generalized  K\"ahler structure}, as a
couple of commuting generalized complex structures $\J_1$ and $\J_2$ such that $\J_1\J_2$ is negative definite;
furthermore, in \cite{Gua-thesis} it is explained the correspondence between generalized K\"ahler structures
and a special type of bihermitian structures which appeared in Theoretical Physics, over twenty years ago \cite{martin}\,.\\
\indent
More precisely, any generalized almost  K\"ahler structure on $M$ corresponds to a quadruple
$(g,b,J_+,J_-)$\,, where $g$ is a Riemannian metric, $b$ is a two-form and $J_{\pm}$ are almost Hermitian
structures on $(M,g)$\,. Furthermore, the corresponding generalized almost K\"ahler structure is integrable if and only if
$J_{\pm}$ are integrable and parallel with respect to $\nabla^{\pm}$, where $\nabla^{\pm}=\nabla\pm\frac 12 g^{-1}h$\,,
with $\nabla$ the Levi-Civita connection of $g$ and $h=\dif\!b$ (equivalently, $J_{\pm}$ are integrable
and $\dif^c_{\pm}\!\o_{\pm}=\mp h$\,, where $\o_{\pm}$ are the K\"ahler forms of $J_{\pm}$).\\
\indent
Classification results for compact bihermitian manifolds were given, mainly in dimension $4$\,, in several papers
(see, for example, \cite{agg}, \cite{ag}\,).\\
\indent
In higher dimensions, a natural case to consider is when $J_+$ and $J_-$ are admissible for an almost quaternionic structure.
This condition was, essentially, considered by physicists who have shown that it holds if and only if the bihermitian structure
is part of a hyperk\"ahler one \cite[Theorem 1]{LRvUZ} (see Theorem~\ref{thm:H^a}\,, below).\\
\indent
By combining this fact with results of \cite{PanWoo-d} and \cite{OrnPan-hologen}\,, we study the `eigendistributions'
of the operator $J_+J_-+J_-J_+$\,. Thus, we obtain natural integrability conditions under which $\dif\!b=0$
(Theorem~\ref{thm:tamed_H^a_integrable}\,, Corollary~\ref{cor:sum_two_H^a_int}\,).

\section{The almost quaternionic generalized K\"ahler manifolds are hyperk\"ahler} \label{section:q_gK}

\indent
A \emph{bundle of associative algebras} is a vector bundle whose typical fibre is an associative algebra $\mathcal{A}$ and whose
structural group is the group of automorphisms of $\mathcal{A}$\,.\\
\indent
An \emph{almost quaternionic structure} on $M$ is a morphism of bundles of associative algebras
$\s:A\to{\rm End}(TM)$\,, where the typical fibre of $A$ is $\Hq$.
Then, $\s({\rm Im}A)$ is an oriented Riemannian vector bundle of rank $3$ and the (local) sections of its sphere bundle
are the \emph{admissible almost complex structures} of $\s$ (see \cite{IMOP}\,).\\
\indent
The following result reformulates \cite[Theorem 1]{LRvUZ}\,. For the reader's convenience, we supply a proof.

\begin{thm} \label{thm:H^a}
Let $(M,L_1,L_2)$ be a generalized almost K\"ahler manifold of dimension at least eight and let $(g,b,J_+,J_-)$
be the corresponding almost bihermititan structure. Suppose that $J_+$ and $J_-$ are admissible almost complex structures
of an almost quaternionic structure on $M$.\\
\indent
Then the following assertions are equivalent:\\
\indent
\quad{\rm (i)} $(M,L_1,L_2)$ is generalized K\"ahler.\\
\indent
\quad{\rm (ii)} $(M,g,J_{\pm})$ are K\"ahler manifolds.\\
\indent
Furthermore, if\/ {\rm (i)} or {\rm (ii)} holds and $J_+\neq\pm J_-$ then the almost quaternionic structure is hyperk\"ahler,
with respect to $(M,g)$\,.
\end{thm}
\begin{proof}
As (ii)$\Longrightarrow$(i) is trivial, it is sufficient to prove that (i)$\Longrightarrow$(ii)\,.\\
\indent
By hypothesis, there exists $a:M\to[-1,1]$ such that $J_+J_-+J_-J_+=-2a$\/ on $M$. If $J_+=\pm J_-$ there is nothing to be proved.
Hence, we may suppose that $a^{-1}\bigl((-1,1)\bigr)\neq\emptyset$\,.\\
\indent
Moreover, as we have to prove that $(M,g,J_{\pm})$ are K\"ahler and, consequently, $a$ is constant, we may assume $a(M)\subseteq(-1,1)$\,.\\
\indent
Then $L_1=L(T^{\C\!}M,\ep_+)$ and $L_2=L(T^{\C\!}M,\ep_-)$\,, where $\ep_{\pm}$ are closed complex two-forms on $M$.
{}From \cite[(6.4) and (6.5)]{Gua-thesis}\,, it quickly follows that
\begin{equation} \label{e:ep_pm}
\begin{split}
({\rm Im}\,\ep_{\pm})(J_+\mp J_-)&=2g\;,\\
({\rm Re}\,\ep_{\pm})(J_+\mp J_-)&=b(J_+\mp J_-)+g(J_+\pm J_-)\;.
\end{split}
\end{equation}
\indent
On multiplying, to the right, both relations of \eqref{e:ep_pm} by $J_+\mp J_-$ we obtain
\begin{equation*}
\begin{split}
(-2\pm 2a)({\rm Im}\,\ep_{\pm})&=2g(J_+\mp J_-)\;,\\
(-2\pm 2a)({\rm Re}\,\ep_{\pm})&=(-2\pm 2a)b\mp g(J_+J_--J_-J_+)
\end{split}
\end{equation*}
and, consequently, $(a-1){\rm Re}\,\ep_+-(a+1){\rm Re}\,\ep_-=-2b$\,.\\
\indent
Therefore
\begin{equation} \label{e:dif_g(J+pmJ-)}
\dif\!\left[\frac{1}{1\pm a}g(J_+\pm J_-)\right]=0\;.
\end{equation}
\indent
Also, as, up to a $B$-field transformation, we may suppose ${\rm Re}\,\ep_-=0$\,, we deduce that the two-form
$\frac{1}{a-1}b$ is closed; equivalently,
\begin{equation} \label{e:dif_b}
\dif\!b=\frac{1}{a-1}\dif\!a\wedge b\;.
\end{equation}
\indent
Note that, the condition $\nabla^{\pm\!}J_{\pm}=0$ is equivalent to
\begin{equation} \label{e:nablaJpm}
g\bigl((\nabla_XJ_{\pm})(Y),Z\bigr)=\mp\tfrac12\bigl[(\dif\!b)(X,J_{\pm}Y,Z)+(\dif\!b)(X,Y,J_{\pm}Z)\bigr]\;,
\end{equation}
for any $X,Y,Z\in TM$.\\
\indent
{}From \eqref{e:dif_b} and \eqref{e:nablaJpm} we obtain
\begin{equation} \label{e:nablaJpm_H^a}
\begin{split}
&g\bigl((\nabla_XJ_{\pm})(Y),Z\bigr)=\pm\frac{1}{2(1-a)}\bigl(\dif\!a\wedge b\bigr)\bigl(X\wedge J_{\pm}Y\wedge Z+X\wedge Y\wedge J_{\pm}Z\bigr)\;,
\end{split}
\end{equation}
for any $X,Y,Z\in TM$.\\
\indent
Obviously, $$K_{\pm}=\frac{1}{\sqrt{2(1\pm a)}}\bigl(J_+\pm J_-\bigr)\;.$$
are anti-commuting almost Hermitian structures on $(M,g)$\,. Furthermore, \eqref{e:nablaJpm_H^a} gives
\begin{equation} \label{e:nablaKpm}
\begin{split}
g\bigl((\nabla_X&K_{\pm})(Y),Z\bigr)=\mp\frac{1}{2(1\pm a)}\,X(a)\,g(K_{\pm}Y,Z)\\
&+\frac{1}{2(1-a)}\left(\frac{1-a}{1+a}\right)^{\pm\frac12}\!\bigl(\dif\!a\wedge b\bigr)\bigl(X\wedge K_{\mp}Y\wedge Z+X\wedge Y\wedge K_{\mp}Z\bigr)\;,
\end{split}
\end{equation}
for any $X,Y,Z\in TM$.\\
\indent
On the other hand, by \eqref{e:dif_g(J+pmJ-)}\,, the almost Hermitian manifolds $\bigl(M,e^{2f_{\pm}}g,K_{\pm}\bigr)$ are $(1,2)$-symplectic,
where $f_{\pm}=-\frac14\log2(1\pm a)$\,. A straightforward calculation shows that this is equivalent to
\begin{equation} \label{e:Kpm_cosymplectic}
\begin{split}
g\bigl((\nabla_{K_{\pm}X}K_{\pm})(Y),Z\bigr)-g\bigl((\nabla&_{X}K_{\pm})(Y),K_{\pm}Z\bigr)=\\
\pm\frac{1}{2(1\pm a)}\bigl[(K_{\pm}Y)(&a)\,g(K_{\pm}X,Z)-(K_{\pm}Z)(a)\,g(K_{\pm}X,Y)\\
&+Y(a)\,g(X,Z)-Z(a)\,g(X,Y)\bigr]\;,
\end{split}
\end{equation}
for any $X,Y,Z\in TM$.\\
\indent
Now, \eqref{e:nablaKpm} and \eqref{e:Kpm_cosymplectic} imply
\begin{equation} \label{e:a_b}
\begin{split}
(K_{\pm}X)(a)\,g(K_{\pm}Y,Z)&+(K_{\pm}Y)(a)\,g(K_{\pm}X,Z)-(K_{\pm}Z)(a)\,g(K_{\pm}X,Y)\\
-X(a)\,&g(Y,Z)+Y(a)\,g(X,Z)-Z(a)\,g(X,Y)=\\
\pm\left(\frac{1-a}{1+a}\right)^{-\frac12}&\bigl(\dif\!a\wedge b\bigr)\bigl(K_{\pm}X\wedge K_{\mp}Y\wedge Z+K_{\pm}X\wedge Y\wedge K_{\mp}Z\\
&-X\wedge K_{\mp}Y\wedge K_{\pm}Z-X\wedge Y\wedge K_{\mp}K_{\pm}Z\bigr)\;,
\end{split}
\end{equation}
for any $X,Y,Z\in TM$.\\
\indent
In \eqref{e:a_b}\,, if from the first relation we subtract the second one, with the roles of $X$ and $Y$ interchanged,
then we obtain
\begin{equation} \label{e:subtracted_a_b}
\begin{split}
&(K_+X)(a)\,g(K_+Y,Z)+(K_+Y)(a)\,g(K_+X,Z)-(K_+Z)(a)\,g(K_+X,Y)\\
&+(K_-X)(a)\,g(K_-Y,Z)+(K_-Y)(a)\,g(K_-X,Z)+(K_-Z)(a)\,g(K_-X,Y)\\
&-2Z(a)\,g(X,Y)=2\left(\frac{1-a}{1+a}\right)^{-\frac12}\!\bigl(\dif\!a\wedge b\bigr)\bigl(K_+X\wedge K_-Y\wedge Z)\;,
\end{split}
\end{equation}
for any $X,Y,Z\in TM$.\\
\indent
{}From \eqref{e:subtracted_a_b}\,, with $Z=K_+X$, it quickly follows that $\grad_{g\!}a$ is zero on the orthogonal complement
of each quaternionic line. As $\dim M\geq8$\,, we obtain that $a$ is constant.
Together with \eqref{e:nablaKpm}\,, this gives that $K_{\pm}$ generate a hyperk\"ahler structure
on $(M,g)$\,, whilst, together with \eqref{e:dif_b}\,, this implies $\dif\!b=0$\,. The proof is complete.
\end{proof}

\begin{rem}
In dimension four, the hypothesis of Theorem \ref{thm:H^a} is equivalent to the condition that $J_+$ and $J_-$ induce the same orientation on $M$,
whilst if $J_+$ and $J_-$ induce different orientations on $M$ then, up to a unique $B$-field transformation, $M$ is locally given by a product
of two K\"ahler manifolds (consequence of \cite[Corollary~5.7]{OrnPan-hologen}\,). Furthermore, there exist four-dimensional generalized K\"ahler
manifolds with $J_+$ and $J_-$ inducing the same orientation and which are not given by a hyperk\"ahler structure (see \cite{Hit-gc_CMP}\,).
\end{rem}

\indent
The next result follows quickly from \eqref{e:dif_b} and \eqref{e:subtracted_a_b}\,.

\begin{cor} \label{cor:gK_dim4}
Let $(M,L_1,L_2)$ be a four-dimensional generalized K\"ahler manifold with $J_+\,,\,J_-$ inducing
the same orientation on $M$ and linearly independent, at each point.\\
\indent
Then, up to a unique $B$-field transformation, the following relations hold:
\begin{equation} \label{e:gK_dim4}
\begin{split}
\dif\!b=&-\frac{1}{1-a}\dif\!a\wedge b\;.\\
*(\dif\!a\wedge b)=&\,\frac{1}{2(1+a)}\,[J_+,J_-](\dif\!a)\;,
\end{split}
\end{equation}
where $*$ is the Hodge star operator of $(M,g)$ and the function $a:M\to(-1,1)$ is characterised by
$J_+J_-+J_-J_+=-2a$\,.
\end{cor}

\indent
We end this section by showing how equations \eqref{e:gK_dim4} can be slightly simplified.

\begin{rem}
Let $(M,L_1,L_2)$ be a four-dimensional generalized K\"ahler manifold with $J_+\,,\,J_-$ inducing
the same orientation on $M$ and linearly independent, at each point.\\
\indent
With the same notations as in Theorem \ref{thm:H^a}, let $K=K_+K_-$\,, $k=\left(\frac{1+a}{1-a}\right)^{\frac 12\!}g$
and $u=\log(1-a)$.\\
\indent
Then \eqref{e:gK_dim4} is equivalent to
\begin{equation} \label{e:gK_dim4_simpl}
\dif\!b=\dif\!u\wedge b=-*_k\!K\!\dif\!u\;.
\end{equation}
\indent
If $\dif\!u$ is nowhere zero, then the second equality of \eqref{e:gK_dim4_simpl} is equivalent to
$$b=c v_{\E}+v_{\F}+\dif\!u\wedge\a\;,$$
where $c$ is a function, $\E$ is generated by $\{\grad u,K(\grad u)\}$, $\F=\E^{\perp}$,
$\a$ is a section of $\F^*$, and $v_{\E}$, $v_{\F}$ are the volume forms of $\E$, $\F$, respectively.
\end{rem}

\section{Factorisation results for generalized K\"ahler manifolds} \label{section:factor_gK}

\indent
Let $(M,L_1,L_2)$ be a generalized K\"ahler manifold and let $(g,b,J_+,J_-)$ be the corresponding bihermitian structure.
For any $a\in[-1,1]$\,, we (pointwisely) denote by $\H^a$ the eigenspace of $J_+J_-+J_-J_+$ corresponding to $-2a$\,; also,
we denote $\H^{\pm}=\H^{\pm1}$ and $\V=\bigl(\H^+\oplus\H^-\bigr)^{\perp}$. Then, at each point of $M$, we have that $\H^a$
are preserved by $J_{\pm}$ and there exist (finite) orthogonal decompositions $TM=\bigoplus_a\H^a$ and $\V=\bigoplus_{|a|<1}\H^a$.

\begin{cor} \label{cor:H^a}
Let $N$ be a complex submanifold of $(M,J_{\pm})$\,, of complex dimension at least four,
endowed with a function $a:N\to(-1,1)$ such that $T_xN\subseteq\H^{a(x)}_{\,x}$\,,
$(x\in N)$\,.\\
\indent
Then $a$ is constant and $N$ is endowed with a natural hyperk\"ahler structure whose underlying Riemannian metric is $g|_N$ and
for which $J_+|_N$ and $J_-|_N$ are admissible complex structures.
\end{cor}
\begin{proof}
As, obviously, $(g,b,J_+,J_-)$ induces a generalized K\"ahler structure on $N$, this follows quickly from Theorem \ref{thm:H^a}\,.
\end{proof}

\indent
{}From \cite[Lemma 2.3]{PanWoo-d} it follows that in an open neighbourhood $U$ of each point of a dense open subset of $M$ there
exist (smooth) functions $a_j:M\to[-1,1]$\,, $(j=1,\ldots,r)$\,, such that $\H^{a_j}$ are distributions on $U$ and $TM=\bigoplus_j\H^{a_j}$;
we call the $\H^{a_j}$ the \emph{(local) eigendistributions} of $J_+J_-+J_-J_+$\,.
Furthermore, if $a$ is a function on $U$ such that, at each point, $-2a$ is an eigenvalue of $J_+J_-+J_-J_+$ then there exists an open subset
of $U$ on which $a=a_j$\,, for some $j$\,; thus, if we assume real-analyticity then $a=a_j$ on $U$.\\
\indent
We point out the following facts:\\
\indent
\quad$\bullet$ \emph{The functions $a_j$ are constant along the integrable manifolds, of dimensions at least eight,
of $\H^{a_j}$, $(j=1,\ldots,r)$\,}; this is a consequence of Corollary \ref{cor:H^a}\,.\\
\indent
\quad$\bullet$ \emph{If $J_+\pm J_-$ are invertible then the holomorphic diffeomorphisms of $(M,L_1,L_2)$ preserve each $\H^{a_j}$,
$(j=1,\ldots,r)$}\,; this is a consequence of \cite[Corollary 6.7]{OrnPan-hologen}\,.

\begin{rem} \label{rem:db=0}
Let $(M,L_1,L_2)$ be a generalized K\"ahler manifold with $\dif\!b=0$\,. Then $(M,g,J_{\pm})$ are K\"ahler and
there exists a nonempty finite subset $A$ of $[-1,1]$ such that, for any $a\in A$\,, we have that $\H^a$ is a parallel foliation
which is holomorphic with respect to both $J_+$ and $J_-$\,. Therefore $(g,J_{\pm})$ induce K\"ahler structures on the leaves of $\H^a$ and,
if $a\neq\pm1$\,, these are admissible with respect to natural hyperk\"ahler structures.
Furthermore, there exist orthogonal decompositions $TM=\bigoplus_{a\in A}\H^a$ and $\V=\bigoplus_{a\in A\setminus\{\pm1\}}\H^a$.\\
\indent
If the cardinal of $A\setminus\{\pm1\}$ is at least two then the leaves of $\bigoplus_{a\in A\setminus\{\pm1\}}\H^a$
are naturally endowed with two distinct hyperk\"ahler structures with respect to which $J_+$ and $J_-$ define admissible complex structures,
respectively.\\
\indent
Furthermore, if $J_++J_-$ (or $J_+-J_-$) is invertible then as, locally, $M$ is the product of a K\"ahler manifold and hyperk\"ahler manifolds,
its holomorphic Poisson structure is the pull-back of the product of the holomorphic symplectic structures of the hyperk\"ahler factors.
\end{rem}

\indent
Next, we prove the following.

\begin{thm} \label{thm:tamed_H^a_integrable}
Let $(M,L_1,L_2)$ be a generalized K\"ahler manifold with $J_++J_-$ (or $J_+-J_-$) invertible
and for which the eigendistributions of $(J_+J_-+J_-J_+)|_{(\H^+\oplus\H^-)^{\perp}}$ have dimensions at least eight.\\
\indent
Then the following assertions are equivalent:\\
\indent
\quad{\rm (i)} $\dif\!b=0$\,.\\
\indent
\quad{\rm (ii)} The eigendistributions of $J_+J_-+J_-J_+$ and their orthogonal complements are integrable.
\end{thm}
\begin{proof}
The implication (i)$\Longrightarrow$(ii) is an immediate consequence of Remark \ref{rem:db=0}\,.\\
\indent
Assume that (ii) holds. {}From \cite[Corollary 6.3]{OrnPan-hologen} it follows that we may suppose that, also, $J_+-J_-$ is invertible.\\
\indent
Then, locally, outside a set with empty interior there exists a finite set $A$ of functions $a:M\to(-1,1)$ such that $\H^a$ are distributions
and $TM=\bigoplus_{a\in A}\H^a$.\\
\indent
Also, $L_1=L(T^{\C\!}M,\ep_+)$ and $L_2=L(T^{\C\!}M,\ep_-)$\,, where $\ep_{\pm}$ are closed complex two-forms on $M$.\\
\indent
By Theorem \ref{thm:H^a}\,, we have that (i) holds if and only if $\dif\!b(X,Y,Z)=0$\,, for any $X\in\H^a$ and
$Y,Z\in\bigoplus_{a'\in A\setminus\{a\}}\H^{a'}$, $(a\in A)$\,.\\
\indent
As $\H^a$, $(a\in A)$\,, are invariant under $B$-field transformations, we may assume ${\rm Re}\,\ep_-=0$\,; equivalently,
$b=-g(J_+-J_-)(J_++J_-)^{-1}$\,.
Together with the fact that $\H^a$, $(a\in A)$\,, and their orthogonal complements are holomorphic foliations, with respect to $J_+$ and $J_-$\,,
this gives that (i) holds if and only if $\H^a$ are Riemannian foliations, $(a\in A)$\,.\\
\indent
Now, note that we, also, have
$${\rm Re}\,\ep_+=b+g(J_++J_-)(J_+-J_-)^{-1}=g\bigl[(J_++J_-)(J_+-J_-)^{-1}-(J_+-J_-)(J_++J_-)^{-1}\bigr]\;.$$
As $L_1$ is integrable, ${\rm Re}\,\ep_+$ is closed and, consequently, $\H^a$ are Riemannian foliations, $(a\in A)$\,.\\
\indent
The proof is complete.
\end{proof}

\indent
We end with the following result.

\begin{cor} \label{cor:sum_two_H^a_int}
Let $(M,L_1,L_2)$ be a generalized K\"ahler manifold for which
the eigendistributions of $(J_+J_-+J_-J_+)|_{(\H^+\oplus\H^-)^{\perp}}$ have dimensions at least eight.\\
\indent
Then the following assertions are equivalent:\\
\indent
\quad{\rm (i)} $\dif\!b=0$\,.\\
\indent
\quad{\rm (ii)} $\H^{\pm}$ and the sum of any two eigendistributions of $J_+J_-+J_-J_+$ are integrable.
\end{cor}
\begin{proof}
The implication (i)$\Longrightarrow$(ii) is trivial.\\
\indent
If (ii) holds then $\H^+\oplus\H^-$ is integrable. Hence, by \cite[Theorem 6.10]{OrnPan-hologen}\,, we may assume $\H^+=0=\H^-$.
The proof follows from Theorem \ref{thm:tamed_H^a_integrable}\,.
\end{proof}


\begin{thebibliography}{10}

\bibitem{agg}
V.~Apostolov, P.~Gauduchon, G.~Grantcharov, Bihermitian surfaces on complex surfaces,
\emph{Proc. London Math. Soc.}, {\bf 79} (1999), 414--428; \emph{Corrigendum}, {\bf 92} (2006), 200--202.

\bibitem{ag}
V.~Apostolov, M.~Gualtieri, Generalized K\"ahler manifolds, commuting complex structures, and split tangent bundle,
\emph{Commun. Math. Phys.}, {\bf 271} (2007), 561--575.

\bibitem{martin}
S.~J.~Gates, C.~M.~Hull, M.~Ro\v{c}ek, Twisted multiplets and new supersymmetric nonlinear sigma models,
\emph{Nuc. Phys. B}, {\bf 248} (1984), 157--186.

\bibitem{Gua-thesis}
M.~Gualtieri, Generalized complex geometry, D.\ Phil.\ Thesis, University of Oxford, 2003.

\bibitem{Hit-gc_CMP}
N.~J.~Hitchin, Instantons, Poisson structures and Generalized K\"ahler Geometry,
\textit{Comm. Math. Phys.}, {\bf 265} (2006) 131-�164.

\bibitem{IMOP}
S.~Ianus, S.~Marchiafava, L.~Ornea, R.~Pantilie, Twistorial maps between quaternionic manifolds,
\textit{Ann. Sc. Norm. Super. Pisa Cl. Sci. (5)}, (to appear).

\bibitem{LRvUZ}
U.~Lindstrom, M.~Ro\v cek, R.~von~Unge, M.~Zabzine, Generalized K\"ahler manifolds and off-shell supersymmetry,
\textit{Commun. Math. Phys.}, {\bf 269} (2007) 833--849.

\bibitem{OrnPan-hologen}
L.~Ornea, R~Pantilie, Holomorphic maps between generalized complex manifolds, Preprint I.M.A.R., Bucharest, 2008,
(\href{http://arxiv.org/abs/0810.1865}{arXiv:0810.1865}).

\bibitem{PanWoo-d}
R.~Pantilie, J.~C.~Wood, Harmonic morphisms with one-dimensional fibres on Einstein manifolds,
\textit{Trans. Amer. Math. Soc.}, {\bf 354} (2002) 4229--4243.



\end{thebibliography}
\end{document}